\documentclass{gtart}

\def\ifplaintex{\expandafter\ifx\csname documentclass\endcsname\relax}

\def\gtp{{\mathsurround=0pt\it $\cal G\mskip-2mu$eometry \&\ 
$\cal T\!\!$opology $\cal P\!$ublications}}  

\def\recd{{\small Received:\qua\receiveddate\ifx\reviseddate\relax
\else\qquad Revised:\qua\reviseddate\fi\par}} 

\def\lognumber#1{\def\thelognumber{#1}}
\def\volumenumber#1{\def\thevolumenumber{#1}}
\def\volumeyear#1{\def\thevolumeyear{#1}}
\def\papernumber#1{\def\thepapernumber{#1}}
\def\pagenumbers#1#2{\def\startpage{#1}\def\finishpage{#2}}
\def\published#1{\def\publishdate{#1}}

\def\received#1{\def\receiveddate{#1}}

\def\accepted#1{\def\accepteddate{#1}}

\long\def\asciiabstract#1{\long\def\theasciiabstract{#1}}

\let\\\par\let\thelognumber\relax\let\thevolumenumber\relax
\let\thepapernumber\relax\let\thevolumeyear\relax\let\startpage\relax
\let\finishpage\relax\let\publishdate\relax\let\receiveddate\relax
\let\reviseddate\relax\let\accepteddate\relax\let\theasciititle\relax
\let\theasciiauthors\relax
\let\theasciiabstract\relax

\let\theasciiemail\relax

\ifplaintex
\font\logobig=cmssbx10 scaled 3836
\font\logomed=cmssbx10 scaled 2557
\else
\font\logobig=cmssbx10 scaled 4200
\font\logomed=cmssbx10 scaled 2800
\fi

\long\def\makeagttitle{   
\count0=\startpage
\agt\hfill      
\hbox to 45truept{\vbox to 0pt{\vglue -13truept{\logomed A\kern -.37em{\logobig 
T}\kern -.38em G}\vss}\hss}
\break
{\small Volume \thevolumenumber\ (\thevolumeyear)
\startpage--\finishpage\nl
Published: \publishdate}

\vglue .25truein

{\parskip=0pt\leftskip 0pt plus
1fil\def\\{\par\smallskip}{\Large\bf\thetitle}\par\medskip} \vglue
0.05truein

{\parskip=0pt\leftskip 0pt plus 1fil\def\\{\par}{\sc\theauthors}
\par\medskip}%
 
\vglue 0.03truein 

{\small\leftskip 25truept\rightskip 25truept{\bf Abstract}\stdspace\theabstract

{\bf AMS Classification}\stdspace\theprimaryclass
\ifx\thesecondaryclass\relax\else; \thesecondaryclass\fi\par
{\bf Keywords}\stdspace \thekeywords\par}\vglue 7truept

}   

\ifplaintex
\font\phead=cmsl9 scaled 950
\font\pnum=cmbx10 scaled 913
\font\pfoot=cmsl9 scaled 950
\headline{\vbox to 0pt{\vskip -4.5mm\line{\small\phead\ifnum
\count0=\startpage ISSN 1472-2739 (on-line) 1472-2747 (printed)
\hfill {\pnum\folio}\else\ifodd\count0\def\\{ }%
\ifx\theshorttitle\relax\thetitle\else\theshorttitle\fi\hfill{\pnum\folio}
\else\def\\{ and }{\pnum\folio}\hfill\ifx\theshortauthors\relax\theauthors
\else\theshortauthors\fi\fi\fi}\vss}}
\footline{\vbox to 0pt{\vglue 0mm\line{\small\pfoot\ifnum\count0=\startpage
\copyright\ \gtp\hfill\else
\agt, Volume \thevolumenumber\ (\thevolumeyear)\hfill\fi}\vss}}
\else
\font=cmsl9 scaled 1050
\font\lnum=cmbx10 
\font=cmsl9 scaled 1050
\makeatletter
\def\@oddhead{{\small\lhead\ifnum\count0=\startpage ISSN 1472-2739 
(on-line) 1472-2747 (printed)\hfill {\lnum\number\count0}\else\ifodd\count0
\def\\{ }\ifx\theshorttitle\relax \thetitle \else\theshorttitle\fi\hfill
{\lnum\number\count0}\else\def\\{ and }{\lnum\number\count0}
\hfill\ifx\theshortauthors\relax 
\theauthors\else\theshortauthors\fi\fi\fi}}\def\@evenhead{\@oddhead}
\def\@oddfoot{\small\lfoot\ifnum\count0=\startpage\copyright\ \gtp\hfill\else
\agt, Volume \thevolumenumber\ (\thevolumeyear)\hfill\fi}
\def\@evenfoot{\@oddfoot}
\makeatother
\fi
\let\maketitlepage\makeagttitle
\let\makeshorttitle\maketitlepage
\let\maketitle\maketitlepage

\newwrite\gtoutfile
\long\gdef\makeheadfile{  
{\def\\{, }\def\s{ }
\immediate\openout\gtoutfile head.xxx
\immediate\write\gtoutfile{To: math@arxiv.org}
\immediate\write\gtoutfile{Subject: put OR rep NNNNN:ppppp}
\immediate\write\gtoutfile{--text follows this line--}
\immediate\write\gtoutfile{Proxy-for: \ifx\theasciiauthors\relax
\theauthors\else\theasciiauthors\fi\s<\ifx\theasciiemail\relax\theemail\else\theasciiemail\fi>}
\immediate\write\gtoutfile{\noexpand\\}
\immediate\write\gtoutfile{Authors: \ifx\theasciiauthors\relax
\theauthors\else\theasciiauthors\fi}
{\def\\{ }\immediate\write\gtoutfile{Title: \ifx\theasciititle\relax
\thetitle\else\theasciititle\fi}}
\immediate\write\gtoutfile{Subj-class: GT or SG, GR etc}
\immediate\write\gtoutfile{MSC-class: \theprimaryclass\ifx\thesecondaryclass\relax\else, \thesecondaryclass\fi}
\immediate\write\gtoutfile{Journal-ref: Algebr. Geom. Topol. \thevolumenumber\s
(\thevolumeyear) \startpage-\finishpage}
\immediate\write\gtoutfile{Comments: Published by Algebraic and
Geometric Topology at}
\immediate\write\gtoutfile{\s\s\s  http://www.maths.warwick.ac.uk/agt/AGTVol\thevolumenumber/agt-\thevolumenumber-\thepapernumber.abs.html}
\immediate\write\gtoutfile{\noexpand\\}
\immediate\write\gtoutfile{}
\ifx\theasciiabstract\relax
\immediate\write\gtoutfile{\theabstract}\else
\immediate\write\gtoutfile{\theasciiabstract}\fi
\immediate\write\gtoutfile{}
\immediate\write\gtoutfile{\noexpand\\}
\immediate\write\gtoutfile{}
\immediate\closeout\gtoutfile}}  

\def\maketitlepage{\makeagttitle\makeheadfile}
\let\makeshorttitle\maketitlepage
\let\maketitle\maketitlepage

\def\ifplaintex{\expandafter\ifx\csname documentclass\endcsname\relax}

\def\gtp{{\mathsurround=0pt\it $\cal G\mskip-2mu$eometry \&\ 
$\cal T\!\!$opology $\cal P\!$ublications}}  

\def\recd{{\small Received:\qua\receiveddate\ifx\reviseddate\relax
\else\qquad Revised:\qua\reviseddate\fi\par}} 

\def\lognumber#1{\def\thelognumber{#1}}
\def\volumenumber#1{\def\thevolumenumber{#1}}
\def\volumeyear#1{\def\thevolumeyear{#1}}
\def\papernumber#1{\def\thepapernumber{#1}}
\def\pagenumbers#1#2{\def\startpage{#1}\def\finishpage{#2}}
\def\published#1{\def\publishdate{#1}}

\def\received#1{\def\receiveddate{#1}}

\def\accepted#1{\def\accepteddate{#1}}

\long\def\asciiabstract#1{\long\def\theasciiabstract{#1}}

\let\\\par\let\thelognumber\relax\let\thevolumenumber\relax
\let\thepapernumber\relax\let\thevolumeyear\relax\let\startpage\relax
\let\finishpage\relax\let\publishdate\relax\let\receiveddate\relax
\let\reviseddate\relax\let\accepteddate\relax\let\theasciititle\relax
\let\theasciiauthors\relax
\let\theasciiabstract\relax

\let\theasciiemail\relax

\ifplaintex
\font\logobig=cmssbx10 scaled 3836
\font\logomed=cmssbx10 scaled 2557
\else
\font\logobig=cmssbx10 scaled 4200
\font\logomed=cmssbx10 scaled 2800
\fi

\long\def\makeagttitle{   
\count0=\startpage
\agt\hfill      
\hbox to 45truept{\vbox to 0pt{\vglue -13truept{\logomed A\kern -.37em{\logobig 
T}\kern -.38em G}\vss}\hss}
\break
{\small Volume \thevolumenumber\ (\thevolumeyear)
\startpage--\finishpage\nl
Published: \publishdate}

\vglue .25truein

{\parskip=0pt\leftskip 0pt plus
1fil\def\\{\par\smallskip}{\Large\bf\thetitle}\par\medskip} \vglue
0.05truein

{\parskip=0pt\leftskip 0pt plus 1fil\def\\{\par}{\sc\theauthors}
\par\medskip}%
 
\vglue 0.03truein 

{\small\leftskip 25truept\rightskip 25truept{\bf Abstract}\stdspace\theabstract

{\bf AMS Classification}\stdspace\theprimaryclass
\ifx\thesecondaryclass\relax\else; \thesecondaryclass\fi\par
{\bf Keywords}\stdspace \thekeywords\par}\vglue 7truept

}   

\ifplaintex
\font\phead=cmsl9 scaled 950
\font\pnum=cmbx10 scaled 913
\font\pfoot=cmsl9 scaled 950
\headline{\vbox to 0pt{\vskip -4.5mm\line{\small\phead\ifnum
\count0=\startpage ISSN 1472-2739 (on-line) 1472-2747 (printed)
\hfill {\pnum\folio}\else\ifodd\count0\def\\{ }%
\ifx\theshorttitle\relax\thetitle\else\theshorttitle\fi\hfill{\pnum\folio}
\else\def\\{ and }{\pnum\folio}\hfill\ifx\theshortauthors\relax\theauthors
\else\theshortauthors\fi\fi\fi}\vss}}
\footline{\vbox to 0pt{\vglue 0mm\line{\small\pfoot\ifnum\count0=\startpage
\copyright\ \gtp\hfill\else
\agt, Volume \thevolumenumber\ (\thevolumeyear)\hfill\fi}\vss}}
\else
\font=cmsl9 scaled 1050
\font\lnum=cmbx10 
\font=cmsl9 scaled 1050
\makeatletter
\def\@oddhead{{\small\lhead\ifnum\count0=\startpage ISSN 1472-2739 
(on-line) 1472-2747 (printed)\hfill {\lnum\number\count0}\else\ifodd\count0
\def\\{ }\ifx\theshorttitle\relax \thetitle \else\theshorttitle\fi\hfill
{\lnum\number\count0}\else\def\\{ and }{\lnum\number\count0}
\hfill\ifx\theshortauthors\relax 
\theauthors\else\theshortauthors\fi\fi\fi}}\def\@evenhead{\@oddhead}
\def\@oddfoot{\small\lfoot\ifnum\count0=\startpage\copyright\ \gtp\hfill\else
\agt, Volume \thevolumenumber\ (\thevolumeyear)\hfill\fi}
\def\@evenfoot{\@oddfoot}
\makeatother
\fi
\let\maketitlepage\makeagttitle
\let\makeshorttitle\maketitlepage
\let\maketitle\maketitlepage

\newwrite\gtoutfile
\long\gdef\makeheadfile{  
{\def\\{, }\def\s{ }
\immediate\openout\gtoutfile head.xxx
\immediate\write\gtoutfile{To: math@arxiv.org}
\immediate\write\gtoutfile{Subject: put OR rep NNNNN:ppppp}
\immediate\write\gtoutfile{--text follows this line--}
\immediate\write\gtoutfile{Proxy-for: \ifx\theasciiauthors\relax
\theauthors\else\theasciiauthors\fi\s<\ifx\theasciiemail\relax\theemail\else\theasciiemail\fi>}
\immediate\write\gtoutfile{\noexpand\\}
\immediate\write\gtoutfile{Authors: \ifx\theasciiauthors\relax
\theauthors\else\theasciiauthors\fi}
{\def\\{ }\immediate\write\gtoutfile{Title: \ifx\theasciititle\relax
\thetitle\else\theasciititle\fi}}
\immediate\write\gtoutfile{Subj-class: GT or SG, GR etc}
\immediate\write\gtoutfile{MSC-class: \theprimaryclass\ifx\thesecondaryclass\relax\else, \thesecondaryclass\fi}
\immediate\write\gtoutfile{Journal-ref: Algebr. Geom. Topol. \thevolumenumber\s
(\thevolumeyear) \startpage-\finishpage}
\immediate\write\gtoutfile{Comments: Published by Algebraic and
Geometric Topology at}
\immediate\write\gtoutfile{\s\s\s  http://www.maths.warwick.ac.uk/agt/AGTVol\thevolumenumber/agt-\thevolumenumber-\thepapernumber.abs.html}
\immediate\write\gtoutfile{\noexpand\\}
\immediate\write\gtoutfile{}
\ifx\theasciiabstract\relax
\immediate\write\gtoutfile{\theabstract}\else
\immediate\write\gtoutfile{\theasciiabstract}\fi
\immediate\write\gtoutfile{}
\immediate\write\gtoutfile{\noexpand\\}
\immediate\write\gtoutfile{}
\immediate\closeout\gtoutfile}}  

\def\maketitlepage{\makeagttitle\makeheadfile}
\let\makeshorttitle\maketitlepage
\let\maketitle\maketitlepage

\def\ifplaintex{\expandafter\ifx\csname documentclass\endcsname\relax}

\def\gtp{{\mathsurround=0pt\it $\cal G\mskip-2mu$eometry \&\ 
$\cal T\!\!$opology $\cal P\!$ublications}}  

\def\recd{{\small Received:\qua\receiveddate\ifx\reviseddate\relax
\else\qquad Revised:\qua\reviseddate\fi\par}} 

\def\lognumber#1{\def\thelognumber{#1}}
\def\volumenumber#1{\def\thevolumenumber{#1}}
\def\volumeyear#1{\def\thevolumeyear{#1}}
\def\papernumber#1{\def\thepapernumber{#1}}
\def\pagenumbers#1#2{\def\startpage{#1}\def\finishpage{#2}}
\def\published#1{\def\publishdate{#1}}

\def\received#1{\def\receiveddate{#1}}

\def\accepted#1{\def\accepteddate{#1}}

\long\def\asciiabstract#1{\long\def\theasciiabstract{#1}}

\let\\\par\let\thelognumber\relax\let\thevolumenumber\relax
\let\thepapernumber\relax\let\thevolumeyear\relax\let\startpage\relax
\let\finishpage\relax\let\publishdate\relax\let\receiveddate\relax
\let\reviseddate\relax\let\accepteddate\relax\let\theasciititle\relax
\let\theasciiauthors\relax
\let\theasciiabstract\relax

\let\theasciiemail\relax

\ifplaintex
\font\logobig=cmssbx10 scaled 3836
\font\logomed=cmssbx10 scaled 2557
\else
\font\logobig=cmssbx10 scaled 4200
\font\logomed=cmssbx10 scaled 2800
\fi

\long\def\makeagttitle{   
\count0=\startpage
\agt\hfill      
\hbox to 45truept{\vbox to 0pt{\vglue -13truept{\logomed A\kern -.37em{\logobig 
T}\kern -.38em G}\vss}\hss}
\break
{\small Volume \thevolumenumber\ (\thevolumeyear)
\startpage--\finishpage\nl
Published: \publishdate}

\vglue .25truein

{\parskip=0pt\leftskip 0pt plus
1fil\def\\{\par\smallskip}{\Large\bf\thetitle}\par\medskip} \vglue
0.05truein

{\parskip=0pt\leftskip 0pt plus 1fil\def\\{\par}{\sc\theauthors}
\par\medskip}%
 
\vglue 0.03truein 

{\small\leftskip 25truept\rightskip 25truept{\bf Abstract}\stdspace\theabstract

{\bf AMS Classification}\stdspace\theprimaryclass
\ifx\thesecondaryclass\relax\else; \thesecondaryclass\fi\par
{\bf Keywords}\stdspace \thekeywords\par}\vglue 7truept

}   

\ifplaintex
\font\phead=cmsl9 scaled 950
\font\pnum=cmbx10 scaled 913
\font\pfoot=cmsl9 scaled 950
\headline{\vbox to 0pt{\vskip -4.5mm\line{\small\phead\ifnum
\count0=\startpage ISSN 1472-2739 (on-line) 1472-2747 (printed)
\hfill {\pnum\folio}\else\ifodd\count0\def\\{ }%
\ifx\theshorttitle\relax\thetitle\else\theshorttitle\fi\hfill{\pnum\folio}
\else\def\\{ and }{\pnum\folio}\hfill\ifx\theshortauthors\relax\theauthors
\else\theshortauthors\fi\fi\fi}\vss}}
\footline{\vbox to 0pt{\vglue 0mm\line{\small\pfoot\ifnum\count0=\startpage
\copyright\ \gtp\hfill\else
\agt, Volume \thevolumenumber\ (\thevolumeyear)\hfill\fi}\vss}}
\else
\font=cmsl9 scaled 1050
\font\lnum=cmbx10 
\font=cmsl9 scaled 1050
\makeatletter
\def\@oddhead{{\small\lhead\ifnum\count0=\startpage ISSN 1472-2739 
(on-line) 1472-2747 (printed)\hfill {\lnum\number\count0}\else\ifodd\count0
\def\\{ }\ifx\theshorttitle\relax \thetitle \else\theshorttitle\fi\hfill
{\lnum\number\count0}\else\def\\{ and }{\lnum\number\count0}
\hfill\ifx\theshortauthors\relax 
\theauthors\else\theshortauthors\fi\fi\fi}}\def\@evenhead{\@oddhead}
\def\@oddfoot{\small\lfoot\ifnum\count0=\startpage\copyright\ \gtp\hfill\else
\agt, Volume \thevolumenumber\ (\thevolumeyear)\hfill\fi}
\def\@evenfoot{\@oddfoot}
\makeatother
\fi
\let\maketitlepage\makeagttitle
\let\makeshorttitle\maketitlepage
\let\maketitle\maketitlepage

\newwrite\gtoutfile
\long\gdef\makeheadfile{  
{\def\\{, }\def\s{ }
\immediate\openout\gtoutfile head.xxx
\immediate\write\gtoutfile{To: math@arxiv.org}
\immediate\write\gtoutfile{Subject: put OR rep NNNNN:ppppp}
\immediate\write\gtoutfile{--text follows this line--}
\immediate\write\gtoutfile{Proxy-for: \ifx\theasciiauthors\relax
\theauthors\else\theasciiauthors\fi\s<\ifx\theasciiemail\relax\theemail\else\theasciiemail\fi>}
\immediate\write\gtoutfile{\noexpand\\}
\immediate\write\gtoutfile{Authors: \ifx\theasciiauthors\relax
\theauthors\else\theasciiauthors\fi}
{\def\\{ }\immediate\write\gtoutfile{Title: \ifx\theasciititle\relax
\thetitle\else\theasciititle\fi}}
\immediate\write\gtoutfile{Subj-class: GT or SG, GR etc}
\immediate\write\gtoutfile{MSC-class: \theprimaryclass\ifx\thesecondaryclass\relax\else, \thesecondaryclass\fi}
\immediate\write\gtoutfile{Journal-ref: Algebr. Geom. Topol. \thevolumenumber\s
(\thevolumeyear) \startpage-\finishpage}
\immediate\write\gtoutfile{Comments: Published by Algebraic and
Geometric Topology at}
\immediate\write\gtoutfile{\s\s\s  http://www.maths.warwick.ac.uk/agt/AGTVol\thevolumenumber/agt-\thevolumenumber-\thepapernumber.abs.html}
\immediate\write\gtoutfile{\noexpand\\}
\immediate\write\gtoutfile{}
\ifx\theasciiabstract\relax
\immediate\write\gtoutfile{\theabstract}\else
\immediate\write\gtoutfile{\theasciiabstract}\fi
\immediate\write\gtoutfile{}
\immediate\write\gtoutfile{\noexpand\\}
\immediate\write\gtoutfile{}
\immediate\closeout\gtoutfile}}  

\def\maketitlepage{\makeagttitle\makeheadfile}
\let\makeshorttitle\maketitlepage
\let\maketitle\maketitlepage

\def\ifplaintex{\expandafter\ifx\csname documentclass\endcsname\relax}

\def\gtp{{\mathsurround=0pt\it $\cal G\mskip-2mu$eometry \&\ 
$\cal T\!\!$opology $\cal P\!$ublications}}  

\def\recd{{\small Received:\qua\receiveddate\ifx\reviseddate\relax
\else\qquad Revised:\qua\reviseddate\fi\par}} 

\def\lognumber#1{\def\thelognumber{#1}}
\def\volumenumber#1{\def\thevolumenumber{#1}}
\def\volumeyear#1{\def\thevolumeyear{#1}}
\def\papernumber#1{\def\thepapernumber{#1}}
\def\pagenumbers#1#2{\def\startpage{#1}\def\finishpage{#2}}
\def\published#1{\def\publishdate{#1}}

\def\received#1{\def\receiveddate{#1}}

\def\accepted#1{\def\accepteddate{#1}}

\long\def\asciiabstract#1{\long\def\theasciiabstract{#1}}

\let\\\par\let\thelognumber\relax\let\thevolumenumber\relax
\let\thepapernumber\relax\let\thevolumeyear\relax\let\startpage\relax
\let\finishpage\relax\let\publishdate\relax\let\receiveddate\relax
\let\reviseddate\relax\let\accepteddate\relax\let\theasciititle\relax
\let\theasciiauthors\relax
\let\theasciiabstract\relax

\let\theasciiemail\relax

\ifplaintex
\font\logobig=cmssbx10 scaled 3836
\font\logomed=cmssbx10 scaled 2557
\else
\font\logobig=cmssbx10 scaled 4200
\font\logomed=cmssbx10 scaled 2800
\fi

\long\def\makeagttitle{   
\count0=\startpage
\agt\hfill      
\hbox to 45truept{\vbox to 0pt{\vglue -13truept{\logomed A\kern -.37em{\logobig 
T}\kern -.38em G}\vss}\hss}
\break
{\small Volume \thevolumenumber\ (\thevolumeyear)
\startpage--\finishpage\nl
Published: \publishdate}

\vglue .25truein

{\parskip=0pt\leftskip 0pt plus
1fil\def\\{\par\smallskip}{\Large\bf\thetitle}\par\medskip} \vglue
0.05truein

{\parskip=0pt\leftskip 0pt plus 1fil\def\\{\par}{\sc\theauthors}
\par\medskip}%
 
\vglue 0.03truein 

{\small\leftskip 25truept\rightskip 25truept{\bf Abstract}\stdspace\theabstract

{\bf AMS Classification}\stdspace\theprimaryclass
\ifx\thesecondaryclass\relax\else; \thesecondaryclass\fi\par
{\bf Keywords}\stdspace \thekeywords\par}\vglue 7truept

}   

\ifplaintex
\font\phead=cmsl9 scaled 950
\font\pnum=cmbx10 scaled 913
\font\pfoot=cmsl9 scaled 950
\headline{\vbox to 0pt{\vskip -4.5mm\line{\small\phead\ifnum
\count0=\startpage ISSN 1472-2739 (on-line) 1472-2747 (printed)
\hfill {\pnum\folio}\else\ifodd\count0\def\\{ }%
\ifx\theshorttitle\relax\thetitle\else\theshorttitle\fi\hfill{\pnum\folio}
\else\def\\{ and }{\pnum\folio}\hfill\ifx\theshortauthors\relax\theauthors
\else\theshortauthors\fi\fi\fi}\vss}}
\footline{\vbox to 0pt{\vglue 0mm\line{\small\pfoot\ifnum\count0=\startpage
\copyright\ \gtp\hfill\else
\agt, Volume \thevolumenumber\ (\thevolumeyear)\hfill\fi}\vss}}
\else
\font=cmsl9 scaled 1050
\font\lnum=cmbx10 
\font=cmsl9 scaled 1050
\makeatletter
\def\@oddhead{{\small\lhead\ifnum\count0=\startpage ISSN 1472-2739 
(on-line) 1472-2747 (printed)\hfill {\lnum\number\count0}\else\ifodd\count0
\def\\{ }\ifx\theshorttitle\relax \thetitle \else\theshorttitle\fi\hfill
{\lnum\number\count0}\else\def\\{ and }{\lnum\number\count0}
\hfill\ifx\theshortauthors\relax 
\theauthors\else\theshortauthors\fi\fi\fi}}\def\@evenhead{\@oddhead}
\def\@oddfoot{\small\lfoot\ifnum\count0=\startpage\copyright\ \gtp\hfill\else
\agt, Volume \thevolumenumber\ (\thevolumeyear)\hfill\fi}
\def\@evenfoot{\@oddfoot}
\makeatother
\fi
\let\maketitlepage\makeagttitle
\let\makeshorttitle\maketitlepage
\let\maketitle\maketitlepage

\newwrite\gtoutfile
\long\gdef\makeheadfile{  
{\def\\{, }\def\s{ }
\immediate\openout\gtoutfile head.xxx
\immediate\write\gtoutfile{To: math@arxiv.org}
\immediate\write\gtoutfile{Subject: put OR rep NNNNN:ppppp}
\immediate\write\gtoutfile{--text follows this line--}
\immediate\write\gtoutfile{Proxy-for: \ifx\theasciiauthors\relax
\theauthors\else\theasciiauthors\fi\s<\ifx\theasciiemail\relax\theemail\else\theasciiemail\fi>}
\immediate\write\gtoutfile{\noexpand\\}
\immediate\write\gtoutfile{Authors: \ifx\theasciiauthors\relax
\theauthors\else\theasciiauthors\fi}
{\def\\{ }\immediate\write\gtoutfile{Title: \ifx\theasciititle\relax
\thetitle\else\theasciititle\fi}}
\immediate\write\gtoutfile{Subj-class: GT or SG, GR etc}
\immediate\write\gtoutfile{MSC-class: \theprimaryclass\ifx\thesecondaryclass\relax\else, \thesecondaryclass\fi}
\immediate\write\gtoutfile{Journal-ref: Algebr. Geom. Topol. \thevolumenumber\s
(\thevolumeyear) \startpage-\finishpage}
\immediate\write\gtoutfile{Comments: Published by Algebraic and
Geometric Topology at}
\immediate\write\gtoutfile{\s\s\s  http://www.maths.warwick.ac.uk/agt/AGTVol\thevolumenumber/agt-\thevolumenumber-\thepapernumber.abs.html}
\immediate\write\gtoutfile{\noexpand\\}
\immediate\write\gtoutfile{}
\ifx\theasciiabstract\relax
\immediate\write\gtoutfile{\theabstract}\else
\immediate\write\gtoutfile{\theasciiabstract}\fi
\immediate\write\gtoutfile{}
\immediate\write\gtoutfile{\noexpand\\}
\immediate\write\gtoutfile{}
\immediate\closeout\gtoutfile}}  

\def\maketitlepage{\makeagttitle\makeheadfile}
\let\makeshorttitle\maketitlepage
\let\maketitle\maketitlepage

\lognumber{31}
\volumenumber{2}
\volumeyear{2002}
\papernumber{31}
\published{9 September 2002}
\pagenumbers{743}{755}
\received{17 May 2002}
\accepted{8 July 2002}

\usepackage{amssymb, amsmath}

\newtheorem{thm}{Theorem}[section] 
\newtheorem*{qthm}{Theorem}
\newtheorem*{qqthm}{Theorem \ref{resfin}} 
\newtheorem*{qqqthm}{Theorem \ref{subgpsep}} 
\newtheorem{lem}[thm]{Lemma}
\newtheorem{cor}{Corollary} 
\theoremstyle{definition}
\newtheorem{defn}[thm]{Definition}

\begin{document}

\title{Engulfing in word-hyperbolic groups}
\authors{G. A. Niblo\\B. T. Williams} 

\address{Faculty of Mathematical Studies, University of 
Southampton\\Highfield, Southampton SO17 1BJ, UK}

\email{G.A.Niblo@maths.soton.ac.uk}

\begin{abstract}   
We examine residual properties of word-hyperbolic groups, ad\-apting a
method introduced by Darren Long to study the residual properties of
Kleinian groups.
\end{abstract}
\asciiabstract{
We examine residual properties of word-hyperbolic groups, adapting a
method introduced by Darren Long to study the residual properties of
Kleinian groups.}


\maketitle

\section {Introduction}

A group is said to be residually finite if the intersection of its
finite index subgroups is trivial. Equivalently it is residually
finite if the trivial subgroup is closed in the profinite topology. It
is an open question whether or not word hyperbolic groups are
residually finite. Evidence that they may be comes from the
observation that many familiar groups in this class are linear and
therefore residually finite by an application of Selberg's
lemma. Furthermore there are geometric methods for establishing the
residual finiteness of free groups \cite{MR11:322e}, surface groups
\cite{MR95b:20050} and some reflection groups \cite{MR58:12996} that
may generalise. Nonetheless the general question seems hard to settle,
hindered by the apparent difficulty of establishing that a given group
contains any proper finite index subgroups at all. In
\cite{MR89h:57001} Long hypothesised this difficulty away by assuming
that the groups he studied satisfied an engulfing property:

\begin{defn} A subgroup $H$ in a group $G$ is said to be engulfed if $H$ is contained in  a proper finite index subgroup of $G$. The group $G$ has the \emph{engulfing property} with respect to a class ${\mathcal H}$ of subgroups of $G$ if every subgroup in the class ${\mathcal H}$ is engulfed in $G$.
\end{defn}

As we will later see Long was able to deduce a strengthened form of
residual finiteness for certain Kleinian groups satisfying a
relatively mild engulfing hypothesis. In \cite{MR2001f:20086} Kapovich
and Wise showed that the question of residual finiteness for the class
of word hyperbolic groups could be reduced to a question concerning
engulfing.

\begin{qthm}[Kapovich, Wise]The following are equivalent:
\begin{itemize}
\item[\rm(i)] Every word-hyperbolic group is residually finite.
\item[\rm(ii)] Every word-hyperbolic group has at least one proper
finite index subgroup.
\end{itemize}
\end{qthm}

The second condition is equivalent to the assertion that every word
hyperbolic group engulfs the identity. While the result of Kapovich
and Wise offers the possibility of an attack on the question of
residual finiteness for the class of word-hyperbolic groups, there is
a real possibility that non-residually finite word hyperbolic groups
exist. In this paper we show how to tackle the more restricted
question of whether a given word-hyperbolic group is residually finite
by suitably adapting Long's method to obtain the following:

\begin{qqthm}
Let $G$ be a word-hyperbolic group and suppose that $G$ engulfs every
finitely generated free subgroup with limit set a proper subset of the
boundary of $G$. Then the intersection of all finite index subgroups
of $G$ is finite. If $G$ is torsion free then it is residually
finite.\end{qqthm}

It is hoped that this result may lead to a new attack on the question
of residual finiteness for certain classes of word hyperbolic groups.

Long's principal aim in introducing engulfing was to establish much
stronger residual properties. A subgroup $H$ of a group $G$ is said to
be \emph{separable} in $G$ if it is an interesection of finite index
subgroups (equivalently $H$ is closed in the profinite topology on
$G$). Residual finiteness is equivalent to separability of the trivial
subgroup.

\begin{qthm}[Long]Let $\Gamma$ be the fundamental group of a closed 
hyperbolic 3-manifold. Suppose that $\Gamma$ has the engulfing
property for those finitely generated subgroups $H$ with $\Lambda
(H)<S^2_\infty$.  Then any geometrically finite subgroup of $\Gamma$
has finite index in its profinite closure.\end{qthm}

There has been substantial recent progress in the field:

\begin{itemize}
\item In \cite{MR2001a:20044} Gitik showed how to construct examples
of closed hyperbolic 3-manifolds such that every quasi-convex subgroup
of the fundamental group is closed in the profinite topology.  Gitik
builds the manifolds by a sequence of doubling operations each of
which consists of glueing two copies of a given compact hyperbolic
3-manifold with non-empty boundary along an incompressible subsurface
of the boundary. Gitik showed that, given appropriate constraints on
the glueing, the fundamental group of the doubled manifold has the
property that all of its quasi-convex subgroups are closed in the
profinite topology. Starting with a handlebody (the fundamental group
of which is free and therefore subgroup separable by Hall's theorem,
\cite{MR11:322e}), Gitik constructs sequences of doubling operations
which yield examples of closed hyperbolic 3-manifolds with fundamental
groups satisfying this property.

\item In \cite{Wise-1} Wise showed that every quasi-convex subgroup of
the fundamental group of the Figure 8 knot complement is closed in the
profinite topology using a geometric method which generalises to many
other link complements, and indeed to other examples arising in
geometric group theory. The conclusion is subsumed by the result of
Long and Reid \cite{MR2002c:57029}.

\item Using arithmetic techniques and building on a method suggested
by the paper of Scott \cite{MR58:12996}, Agol, Long and Reid
\cite{MR1836283} showed that the geometrically finite subgroups of
Bianchi groups are closed in the profinite topology.
\end{itemize}

In our second main result we again adapt Long's technique to show:

\begin{qqqthm}Let $G$ be a  word-hyperbolic group which engulfs every 
finitely generated subgroup $K$ such that the limit set $\Lambda(K)$
is a proper subset of the boundary of $G$. Then every quasi-convex
subgroup of $H$ has finite index in its profinite closure in
$G$. \end{qqqthm}

It may be that existing proofs of separability can be simplified using
this result, but by way of caution we also generalise a construction
of Long's to show that every non-elementary word hyperbolic group
contains proper subgroups which fail to be engulfed. However the
construction sheds no light on the question of engulfing for finitely
generated subgroups.

The work of adapting Long's argument to the context of torsion free
word hyperbolic groups formed part of the thesis of the second author
\cite{Williams-1}. The main technical difficulties in this paper arise
in adapting the argument to the presence of torsion.

\section {Word-hyperbolic groups}

This section is a brief introduction to word-hyperbolic groups. The
reader is referred to \cite{MR92f:53050} for a full treatment.

Let $G$ be a finitely generated group, let $S$ be a finite generating
set for $G$, and consider $G$ as a metric space with respect to the
word metric corresponding to this generating set.

The group $G$ is said to be \emph{word-hyperbolic} if it is a
$\delta$-hyperbolic space for some $\delta\geq 0$.

 The \emph{boundary at infinity} of $G$, denoted $\partial G$ is
defined as a metric space whose points are equivalence classes of rays
converging to infinity in the group. It is the dynamics of the action
of $G$ (and its subgroups) on this boundary that we will use to prove
the main theorems in this paper. We take a moment to recall the
important features of the boundary and of those dynamics.

 A word-hyperbolic group is called \emph{elementary} if it is finite
 or contains a finite index infinite cyclic subgroup and is
 \emph{non-elementary} otherwise. An elementary word-hyperbolic group
 is either finite, in which case it has an empty boundary at infinity,
 or it is virtually cyclic in which case its boundary consists of two
 points.  For any word hyperbolic group the boundary is compact and
 metrisable, and non-elementary word-hyperbolic groups have infinite
 boundaries in which there are no isolated points.

Given a subgroup $H$ of $G$, the \emph{limit set of $H$} which is
denoted $\Lambda (H)$ is defined as the subset of $\partial G$
attainable by sequences of elements of $H$. $H$ acts properly
discontinuously on $\partial G \setminus \Lambda (H)$.
 
The following describes the action of infinite order elements on the
boundary. If $g$ is an infinite order element of $G$ it acts on the
Cayley graph $ G$ by translation along a quasi-geodesic line, $\alpha$
say, (obtained by joining $g^i$ to $g^{i+i}$ for all $i\in{\mathbb Z}$
by a geodesic in $G$). Denote by $\partial g=\{\partial g^+,\partial
g^- \}=\{ lim_{i\rightarrow\infty}g^i,
lim_{i\rightarrow\infty}g^{-i}\}$ the endpoints of $\alpha$ in
$\partial G$ (which are fixed by $g$). There exist disjoint
neighbourhoods $U_+$ and $U_-$ of $\partial g^+$ and $\partial g^-$
respectively such that for sufficiently large $r$ and all $x\in
\partial G\setminus (U_+\cup U_-)$ we have $g^rx\in U_+$ and
$g^{-r}x\in U_-$. We say that the pair $(U_+,U_-)$ is \emph{absorbing}
for $g^r$. In fact any pair of disjoint neighbourhoods of $\partial
g^+$ and $\partial g^-$ is absorbing for $g^k$ for sufficiently large
$k$. (See \cite{MR92f:53050} Chapter 8.)

The following well known fact can be viewed as an alternative
definition of the limit set of a subgroup. A proof is included for the
convenience of the reader.

\begin{lem}\label{smallest} 
Let $H$ be a non-elementary subgroup of a word-hyperbolic group
$G$. Then $\Lambda (H)$ is the smallest non-empty closed $H$-invariant
subset of $\partial G$.
\end{lem}

\begin{proof}
We prove that if $A\subseteq\partial G$ is closed and $H$-invariant
then $\Lambda (H)\subseteq A$. Firstly, let $B\subset \partial
G$. Denote by $I(B)$ the set of points of $G$ lying on geodesics
between points of $B$. Suppose that $B\not=\emptyset$ and $|B|\not=
1$. Then $I(B)\not=\emptyset$. Let $\{ x_i\}\subseteq I(B)$ be a
sequence such that $x_i\rightarrow x\in \partial G$. We claim that
$x\in \overline B$. To see this, for each $i$ choose a geodesic
$l_i=[b_i',b_i'']$ with $b_i',b_i''\in B$. Passing to a subsequence if
necessary we get $b_i'\rightarrow b'\in \overline B$,
$b_i''\rightarrow b''\in \overline B$, $l_i\rightarrow l$.
$x_i\rightarrow x\in l\cup \{b',b''\}$ and hence $x\in \{b',b''\}$.

Now let $A\subseteq \partial G$ be closed and $H$-invariant. Let
$I(A)$ be as above. Then $I(A)$ is $H$-invariant. First suppose that
$1\in I(A)$. Then $H\subseteq I(A)$. Let $x\in \Lambda (H)$ and $\{
x_i\}\subseteq H$ so that $x_i\rightarrow x$. By the first paragraph
of the proof $x\in {\overline A}=A$ and hence $\Lambda (H)\subseteq
A$. Now suppose that $1\not\in I(A)$. Then $I(A)\cap H=\emptyset$ and
$I(A)$ is a union of right cosets of $H$.  Suppose that $Hg\subseteq
I(A)$. Let $x\in\Lambda (H)$ and $\{ x_i\}\subseteq H$ with
$x_i\rightarrow x$. Then since $x_ig$ and $x_i$ are a distance exactly
$|g|$ apart for all $i$ we have $x_ig\rightarrow x\in \Lambda (H)$ and
hence $x\in {\overline A}=A$ and $\Lambda (H)\subseteq A$ as required.

It is clear that $\Lambda (H)$ is $H$-invariant so it remains to prove
that $\Lambda (H)$ is closed. We show that $\partial {G }\setminus
\Lambda (H)$ is open. Let $y\in\partial {G }\setminus \Lambda (H)$ and
let $\{ y_i\}$ be a sequence converging to $y$. Let $\alpha_i$ be
geodesics realising the distances $d(y_i,H)$. There is no bound on the
lengths of the $\alpha_i$. Let $z_i$ lie on $\alpha_i$ so that there
is no bound on the distances $d(y_i,z_i)$ and $d(z_i,H)$. Let $\{
z_i\}$ converge to $z\in \partial G$ then the horoball $N_{(y.z)}(y)$
is an open set containing $y$ and disjoint from $\Lambda (H)$ as
required.
\end{proof}

\begin{cor}\label{limset}
 Let $H$ be a non-elementary subgroup of a word-hyperbolic group $G$.
Then $\Lambda (H)$ is the closure of the set
$$ S=\{\partial h^+, \partial h^- | h\in H\hbox{, $h$ has infinite
order}\}\subset\partial G.$$
\end{cor}

\begin{proof} By Lemma \ref{smallest} $\Lambda (H)$ is the minimum 
non-empty closed $H$-invariant subset of $\partial G$. For any
infinite order element $h\in H$, the limit points $\partial h^\pm$
both lie in $\Lambda H$, hence the closure of $S$ must be contained in
$\Lambda (H)$. On the other hand $S$ is clearly $H$-invariant and so
by Lemma \ref{smallest} its closure contains $\Lambda(H)$ as required.
\end{proof}

We will need the following technical observation:

\begin{lem}\label{torfreegens}
Let G be a non-elementary word-hyperbolic group with generators
$g_1,\ldots g_n$, and $N$ a subgroup of $G$ with $\Lambda N=\partial
G$. Then there are infinite order elements $x_1,\ldots x_n$ in $N$
such that the elements $x_ig_ix_i$ generate a free subgroup $H<G$ with
$\Lambda H\not=\partial G$. In particular $G$ is generated by the
subset $\{x_1,\ldots x_n, x_1g_1x_1,\ldots x_ng_nx_n\}$ which consists
of elements of infinite order.
\end{lem}

\begin{proof}
Since $G$ is non-elementary its boundary is infinite, and since limit
sets of infinite order elements of $N$ are dense, given any non-empty
open subset $\mathcal U$ of the boundary we may choose elements
$y_1,\ldots,y_n\in N$ of infinite order with $\partial y_i\subseteq
\mathcal U$ for all $i$ and $\partial y_i\cap \partial y_j=\emptyset$
for $i\not = j$. In particular if we let $\mathcal U$ be the
complement in $\partial G$ of the union of the fixed sets of the
infinite order elements in the set $\{g_1,\ldots, g_n\}$ then we can
also ensure that $\partial g_i\cap \partial y_j=\emptyset$ for all
$i,j$

If a generator $g_i$ acts trivially on the boundary then set
$x_i=y_i$. The element $x_ig_ix_i$ acts on the boundary in the same
way as the infinite order element $y_i^2$, and its two fixed points
are $\partial y_i^\pm\in \mathcal U$. If the generator $g_i$ has
infinite order then since its fixed points are disjoint from those of
$y_i$ (and the boundary is metrisable), we may choose small
neighbourhoods $\mathcal U_i^\pm$ of the limit points $\partial
y_i^\pm$ such that $g_i^{\pm 1}(\mathcal U_i^+\cup \mathcal U_i^-)\cap
(\mathcal U_i^+\cup \mathcal U_i^-)=\emptyset$. We choose the
neighbourhoods $\mathcal U_i^\pm$ small enough to be disjoint and so
that the complement of the closure of the union of the neighbourhoods
is non-empty.

The neighbourhoods $\mathcal U_i^\pm$ are absorbing for any
sufficiently high power $y_i^{\pm r}$ of $y_i$, and it follows easily
that setting $x_i=y_i^r$ the neighbourhoods are absorbing for
$(x_ig_ix_i)^{\pm 1}$.  To see this choose any point $p$ in the
complement of $\mathcal U_i^+\cup \mathcal U_i^-$. Its image $x_i(p)$
lies in $\mathcal U_i^+$, and since $g_i(\mathcal U_i^+\cup \mathcal
U_i^-)\cap (\mathcal U_i^+\cup \mathcal U_i^-)=\emptyset$ $g_ix_i(p)$
does not lie in $\mathcal U_i^+\cup \mathcal U_i^-$. Hence
$x_ig_ix_i(p)$ lies in $\mathcal U_i^+$.  A similar argument shows
that $x_i^{-1}g_i^{-1}x_i^{-1}(p)$ lies in $\mathcal U_i^-$, and
iterating shows that $(x_ig_ix_i)^r(p)\in \mathcal U_i^+\cup \mathcal
U_i^-$ for any non-zero power of the element $x_ig_ix_i$.

We will now use the standard Schottky argument to show that these
elements generate a free subgroup. Let
$w=x_{i_1}^{\epsilon_{i_1}}g_{i_1}^{\epsilon_{i_1}}x_{i_1}^{\epsilon_{i_1}}\ldots
x_{i_s}^{\epsilon_{i_s}}g_{i_s}^{\epsilon_{i_s}}x_{i_s}^{\epsilon_{i_s}}$
be a reduced word in the elements $x_ig_ix_i$ and their inverses, and
choose a point $p$ in the complement of the union of the absorbing
pairs $\mathcal U_i^+\cup \mathcal U_i^-$. As argued above
$x_{i_s}^{\epsilon_{i_s}}g_{i_s}^{\epsilon_{i_s}}x_{i_s}^{\epsilon_{i_s}}(p)\in
\mathcal U_{i_s}^+\cup \mathcal U_{i_s}^-$. If $i_{s-1}=i_s$ then we
may iterate to see that the image of $p$ under the element
$x_{i_{s-1}}^{\epsilon_{i_{s-1}}}g_{i_{s-1}}^{\epsilon_{i_{s-1}}}x_{i_{s-1}}^{\epsilon_{i_{s-1}}}$
also lies in $\mathcal U_{i_s}^+=\mathcal U_{i_{s-1}}^+$. If
$i_{s-1}\not= i_s$ then, since the absorbing set $\mathcal
U_{i_{s}}^+\cup \mathcal U_{i_{s}}^-$ is disjoint from the absorbing
set $\mathcal U_{i_{s-1}}^+\cup \mathcal U_{i_{s-1}}^-$, the image
$x_{i_{s-1}}^{\epsilon_{i_{s-1}}}g_{i_{s-1}}^{\epsilon_{i_{s-1}}}x_{i_{s-1}}^{\epsilon_{i_{s-1}}}\ldots
x_{i_s}^{\epsilon_{i_s}}g_{i_s}^{\epsilon_{i_s}}x_{i_s}^{\epsilon_{i_s}}(p)$
lies in $\mathcal U_{i_{s-1}}^+\cup \mathcal U_{i_{s-1}}^-$. Iterating
the argument we see that the point $p$ ends in the absorbing pair
$\mathcal U_{i_1}^+\cup \mathcal U_{i_1}^-$. Since it did not start
there it is not invariant under the action of the element $w$ which is
therefore not the identity. Hence every reduced word in the generators
$x_ig_ix_i$ is non-trivial and the subgroup is free as
required. Finally we note that since the accumulation points for the
action of this subgroup $H$ lie in the union of the absorbing pairs
$\mathcal U_{i}^+\cup \mathcal U_{i}^-$ the limit set of this subgroup
lies in the closure of their union. Since this closure is not all of
$\partial G$ neither is $\Lambda H$.
\end{proof}

\section{Separability}

The profinite topology on a group $G$ is defined by taking as a basis
for the closed sets the cosets of all finite index normal subgroups of
$G$.  Note that finite index subgroups are themselves closed, and
(since the complement is a finite union of cosets each of which is
also open) they are also open. Given a subgroup $H<G$ we will denote
the closure of $H$ in the profinite topology on $G$ by $\overline H$.

\begin{defn} Given a group $G$, a finitely generated subgroup $H$ is 
\emph{separable in $G$} if it is closed in the profinite topology on
$G$. A group $G$ is \emph{residually finite} if $\{ e\}$ is closed and
is \emph{subgroup separable} or \emph{LERF (locally extended
residually finite)} if every finitely generated subgroup $H$ is
separable in $G$. A word-hyperbolic group is qc subgroup separable if
every quasi convex subgroup is closed in the profinite
topology. \end{defn}

Note that if a group is subgroup separable then a fortiori it has the
engulfing property for its finitely generated subgroups. On the other
hand in \cite{MR99k:57041} examples are given of fundamental groups of
geometric 3-manifolds which contain two generator subgroups which are
not even engulfed. These examples, based on earlier examples of
\cite{MR88g:20057} are not word-hyperbolic, however Long showed in
\cite{MR89h:57001} that the fundamental group of a hyperbolic
3-manifold always contains (infinitely generated) subgroups that are
not engulfed.

\section{(Almost) residual finiteness}

For this section let $N$ denote the residual core of $G$, i.e.,
intersection of all finite index subgroups. (This is the closure
$\overline{\{e\}}$ of the trivial subgroup in the profinite topology.)
This subgroup is normal and therefore \cite{MR92f:53050} its limit set
is either empty (if $N$ is finite) or all of $\partial G$ (if $N$ is
infinite).We will say that the group $G$ is \emph{almost residually
finite} if $N$ is a finite subgroup. Note that torsion free almost
residually finite groups are residually finite.

\begin{thm}\label{resfin} Let $G$ be a word-hyperbolic group and
suppose that $G$ engulfs every finitely generated free subgroup $S$
such that $\Lambda(S)$ is a proper subset of $\partial G$. Then $G$ is
almost residually finite. If $G$ is torsion free, then it is
residually finite.\end{thm}

\begin{proof}
If $G$ is elementary then it is either finite or virtually cyclic. In
both cases it is trivially residually finite, so we may assume that
$G$ is non-elementary and has infinite boundary.

Let $\{g_i\mid 1\leq i\leq n\}$ be a generating set for $G$. If $G$ is
 not almost residually finite then $\Lambda (N)=\partial G$. It
 follows from Lemma \ref{torfreegens} that we may choose elements
 $x_i\in N$ such that the elements $x_ig_ix_i$ generate a free
 subgroup $H$ with $\Lambda H\not=\partial G$. By hypothesis $H$ is
 engulfed, so there is a proper finite index subgroup $L<G$ with
 $H<L$. The subgroup $L$ must contain the elements $x_ig_ix_i$, but by
 hypothesis $N<L$ so it also contains the elements $x_i$. Hence it
 contains all of the generators $g_i$ of $G$. This is a contradiction.
 Hence $G$ is almost residually finite, and if $G$ is torsion free it
 is residually finite.
\end{proof}

\section{(Almost) subgroup separability}

Note that if $H$ is a finite subgroup of an almost residually finite
group $G$, and if $N$ is the intersection of the finite index
subgroups of $G$, then $HN$ is finite, and is closed. Hence the
intersection of the finite index subgroups of $G$ containing $H$ is a
finite extension of $H$.

\begin{defn}
We will say that a subgroup $H<G$ is \emph{almost separable} if $H$
has finite index in $\overline H$.
\end{defn}

\begin{thm}\label{subgpsep} Let $G$ be a  non-elementary word-hyperbolic group. Suppose
that $G$ has the engulfing property for all finitely generated
subgroups $K$ such that $\Lambda(K)$ is a proper subset of $\partial
G$. Then every quasi-convex subgroup $H\leq G$ is almost separable in
$G$.
\end{thm}

\begin{proof} 
Applying Theorem \ref{resfin} we see that the intersection $N$ of all
finite index subgroups of $G$ is finite. It is easy to see that $G/N$
is itself residually finite.

Let $KN/N$ be any subgroup of $G/N$ with limit set a proper subset of
the boundary of $G/N$. There is a G-equivariant quasi-isometry from
$G$ to $G/N$ taking $KN$ to $KN/N$ and it follows that the limit set
of $KN$ is a proper subset of the boundary of $G$. By the hypothesis
there is a proper finite index subgroup of $G$ containing $KN$, and
since it contains $N$ its image is a proper finite index subgroup of
$G/N$ containing $KN/N$. Hence $G/N$ satisifes the hypotheses of the
theorem, but in addition it is residually finite.

Now suppose the theorem is true for $G/N$. Let $H$ be a quasi-convex
subgroup of $G$, so $HN/N$ is a quasi-convex subgroup of $G/N$. By the
assumption, $HN/N$ has finite index in its closure $\overline{HN/N}$
under the profinite topology. Since the map $G\longrightarrow G/N$ is
continuous the preimage of $\overline{HN/N}$ is itself closed in $G$
and clearly contains $H$ as a subgroup of finite index. Hence in order
to establish the theorem for $G$ it suffices to establish it for
$G/N$. This reduces us to the case where $G$ is residually finite, so
from now on we make this additional assumption.

Now since $G$ is residually finite, its finite subgroups and its
maximal abelian subgroups (see \cite{MR92g:20047}) are all closed in
the profinite topology. Since $G$ is word-hyperbolic its maximal
abelian subgroups are virtually cyclic, and therefore every elementary
subgroup of $G$ has finite index in its profinite closure. Hence we
can assume that $H$ is non-elementary.

We will make use of the following observation. A proof is given in
Kapovich and Short \cite{MR98g:20060}.

 \begin{lem} Let $H$ be a quasiconvex subgroup of a word-hyperbolic
group $G$. If $H<L<G$ with $\Lambda (H)=\Lambda (L)$ then
$|H:L|<\infty$.
\end{lem}

It follows from this that it suffices to show that the profinite
closure $\overline{H}$ of any non-elementary quasi-convex subgroup
$H<G$ has the same limit set as $H$. For the remainder of the argument
fix a generating set $\{ g_1,g_2,\ldots ,g_n\}$ for $G$. By Lemma
\ref{torfreegens} we can choose this set to consist of infinite order
elements.

Since $H\leq \overline H$ clearly $\Lambda (H)\subseteq \Lambda
(\overline H)$, and if $\Lambda (H)=\partial G$ then the result is
clear so suppose that $\Lambda (H)$ is a proper subset of $\partial
G$. Assume, for a contradiction, that $\Lambda (H)\not=\Lambda
(\overline H)$.

Choose a point $p\in \Lambda (\overline H)\setminus \Lambda (H)$. By
Corollary \ref{limset} there is a sequence of infinite order elements
$k_i\in \overline{H}$ with fixed points
$p_i\in\Lambda(\overline{H})\subset \partial G$ such that the sequence
$p_i$ converges to $p$. Since $\Lambda(H)$ is closed and $p\not\in
\Lambda H$ almost all the points $p_i$ are also not in $\Lambda(H)$,
so almost all the elements $k_i$ are in $\overline{H}\setminus H$ and,
since limit sets of non-elementary quasi-convex subgroups have no
isolated limit points, without loss we can choose them to have
distinct limit sets. Hence we can choose one of them with limit points
$p^\pm$ in $\partial G$ distinct from the limit points of the
generators. Since $p^\pm$ are also not in $\Lambda(H)$ we may choose
an absorbing pair of neighbourhoods $U^\pm$ of the pair $p^\pm$
disjoint from $\Lambda(H)\cup\{ g_1,g_2,\ldots ,g_n\}$. Since $G$ acts
uniformly on its boundary and $\Lambda H$ is a closed set disjoint
from the limit points of $k$, for some power $k^r$ the image
$k^r(\Lambda(H))$ is contained in $U^+$ and is therefore disjoint from
$\Lambda(H)\cup\{ g_1,g_2,\ldots ,g_n\}$. The image $k^r(\Lambda(H))$
is the limit set of $k^r Hk^{-r}$ which by construction is a subgroup
of $\overline{H}$.

Since $H$ is non-elementary so is $k^r Hk^{-r}$ and we may choose
elements $y_1,y_2,$ $\ldots,y_n\in k^r Hk^{-r}$ with distinct fixed
sets in the boundary. Notice that by our construction of the subgroup
$k^r Hk^{-r}$ the fixed points $\partial y_1,\partial
y_2,\ldots,\partial y_n$ lie in $\Lambda (\overline H)-\Lambda (H)$
and $\partial y_i\not= \partial g_j$ for any $i,j$.  We may later need
to modify the choice of these elements by taking powers of them. In
doing so we do not change their fixed points.

Let $C\subset \partial G - \Lambda (H)$ be a compact set containing
the fixed points of the elements $y_i$ in its interior (the closure of
a sufficiently small open metric ball around the fixed points will
do). $H$ acts properly discontinuously on $\partial G-\Lambda (H)$ so
there are finitely many non-trivial elements of $H$, $h_1,h_2,\ldots
h_m$ say, taking $C$ to intersect itself.  Since $G$ is residually
finite so is $H$, and so there exists a finite index normal subgroup
$A\triangleleft H$ containing none of the $h_i$.

We now need the following technical Lemma taken from
\cite{MR89h:57001}.

\begin{lem}\label{powers} Let $G$ and $H$ be as above and suppose that $A\triangleleft H$ is a normal subgroup of  index  $t$ in $H$. For any element  $h\in \overline H$, $h^t\in \overline A$. \end{lem}
 
Since taking powers of the elements $y_i$ does not change their fixed
points we can use this lemma to ensure that the elements $y_i$ all lie
in the subgroup $\overline A$.  Since $\partial G$ is metrisable we
can choose $n$ mutually disjoint pairs of neighbourhoods
$(U_i^+,U_i^-)$ for the $\partial y_i$ so that the closure of each is
contained in the interior of $C$. Ensure that $(U_i^+,U_i^-)$ is
absorbing for $y_i$ by again taking large powers and relabelling.

Now let $s_i=y_ig_iy_i$ and consider the group $B$ generated by the
elements $s_i$ together with the generators of $A$. Since $A$ has
finite index in the finitely generated group $H$ it too is finitely
generated and so is $B$. We claim that its limit set is contained in
the closure of $\cup_i(U_i^+,U_i^-)\cup(\partial G -C)$.

Let $U_i=U_i^+\cup U_i^-$.

The limit set is the closure of the H-orbit of any point in it (by
\ref{limset}). Choose a point $p\in C- \cup_i U_i$ and write an
arbitrary element $b\in B$ as a reduced word
$s_{i_1}^{\epsilon_1}a_1s_{i_2}^{\epsilon_2}a_2\ldots
s_{i_k}^{\epsilon_k}a_k$ where $a_i\in A$, where possibly $s_{i_1}$ or
$a_k$ may be the identity elements, but none of the other elements are
trivial. We examine the image of $p$ under the action of $b$; there
are four cases to consider:

Neither $s_{i_1}$ nor $a_k$ is the identity:

$b(p)= s_{i_1}^{\epsilon_1}a_1s_{i_2}^{\epsilon_2}a_2\ldots
s_{i_k}^{\epsilon_k}a_k(p) \in
s_{i_1}^{\epsilon_1}a_1s_{i_2}^{\epsilon_2}a_2\ldots
s_{i_k}^{\epsilon_k}(\partial G -C)$\nl 
\hbox{}\qquad$\subset
s_{i_1}^{\epsilon_1}a_1s_{i_2}^{\epsilon_2}a_2\ldots
s_{i_k}^{\epsilon_k}(\partial G -(U_{i_k})) \subset
s_{i_1}^{\epsilon_1}a_1s_{i_2}^{\epsilon_2}a_2\ldots
a_{i_{k-1}}U_{i_k}$\nl 
\hbox{}\qquad$\subset
s_{i_1}^{\epsilon_1}a_1s_{i_2}^{\epsilon_2}a_2\ldots a_{i_{k-1}}C
\subset s_{i_1}^{\epsilon_1}a_1s_{i_2}^{\epsilon_2}a_2\ldots
s_{i_{k-1}}^{\epsilon_{k-1}}(\partial G -C)$\nl 
\hbox{}\qquad\qquad$\smash{\vdots}\vrule height 10pt width 0pt$\nl
\hbox{}\qquad$\subset
s_{i_1}^{\epsilon_1}a_1C\subset s_{i_1}^{\epsilon_1}(\partial G-C)
\subset s_{i_1}^{\epsilon_1}(\partial G-(U_{i_1})) \subset U_{i_1}$

Only $a_k$ is the identity:

$b(p)= s_{i_1}^{\epsilon_1}a_1s_{i_2}^{\epsilon_2}a_2\ldots
s_{i_k}^{\epsilon_k}a_k(p) \in
s_{i_1}^{\epsilon_1}a_1s_{i_2}^{\epsilon_2}a_2\ldots
s_{i_k}^{\epsilon_k}(C- \cup_i U_i)$\nl 
\hbox{}\qquad$\subset
s_{i_1}^{\epsilon_1}a_1s_{i_2}^{\epsilon_2}a_2\ldots
a_{i_{k-1}}U_{i_k} \subset
s_{i_1}^{\epsilon_1}a_1s_{i_2}^{\epsilon_2}a_2\ldots a_{i_{k-1}}C$\nl 
\hbox{}\qquad$\subset s_{i_1}^{\epsilon_1}a_1s_{i_2}^{\epsilon_2}a_2\ldots
s_{i_{k-1}}^{\epsilon_{k-1}}(\partial G -C)$\nl
\hbox{}\qquad\qquad$\smash{\vdots}\vrule height 10pt width 0pt$\nl
\hbox{}\qquad$\subset
s_{i_1}^{\epsilon_1}a_1C \subset s_{i_1}^{\epsilon_1}(\partial G-C)
\subset s_{i_1}^{\epsilon_1}(\partial G-(U_{i_1})) \subset U_{i_1} $

Only $s_{i_1}$ is the identity:

$b(p)= a_1s_{i_2}^{\epsilon_2}a_2\ldots s_{i_k}^{\epsilon_k}a_k(p) \in
a_1s_{i_2}^{\epsilon_2}a_2\ldots s_{i_k}^{\epsilon_k}(C- \cup_i U_i)
\subset a_1s_{i_2}^{\epsilon_2}a_2\ldots a_{i_{k-1}}U_{i_k}$\nl
\hbox{}\qquad$\subset
a_1s_{i_2}^{\epsilon_2}a_2\ldots a_{i_{k-1}}C \subset
a_1s_{i_2}^{\epsilon_2}a_2\ldots s_{i_{k-1}}^{\epsilon_{k-1}}(\partial
G -C)$\nl 
\hbox{}\qquad\qquad$\smash{\vdots}\vrule height 10pt width 0pt$\nl
\hbox{}\qquad$\subset
a_1C \subset (\partial G-C) $

Both $s_{i_1}$ and $a_k$ are the identity:

$b(p)= a_1s_{i_2}^{\epsilon_2}a_2\ldots s_{i_k}^{\epsilon_k}a_k(p) \in
a_1s_{i_2}^{\epsilon_2}a_2\ldots s_{i_k}^{\epsilon_k}(C- \cup_i U_i)
\subset a_1s_{i_2}^{\epsilon_2}a_2\ldots a_{i_{k-1}}U_{i_k}$\nl
\hbox{}\qquad$\subset
a_1s_{i_2}^{\epsilon_2}a_2\ldots a_{i_{k-1}}C \subset
a_1s_{i_2}^{\epsilon_2}a_2\ldots s_{i_{k-1}}^{\epsilon_{k-1}}(\partial
G -C)$\nl 
\hbox{}\qquad\qquad$\smash{\vdots}\vrule height 10pt width 0pt$\nl
\hbox{}\qquad$\subset
a_1C \subset (\partial G-C) $

The conclusion is that $p$ ends up in $\cup_i U_i$ or in $\partial
G-C$, and in particular the closure of its orbit lies in the union of
the closures of these subsets as required.

Hence $B$ is a finitely generated subgroup of $G$ with $\Lambda(B)$ a
proper subset of $\partial G$ and our engulfing hypothesis for such
subgroups ensures that there exists a proper finite index subgroup
$K<G$ containing $B$. Since this subgroup contains $A$ it also
contains $\overline A\leq K$ and hence $K$ contains the elements
$y_1,y_2,\ldots ,y_n$. But $K$ also contains the elements
$s_i=y_ig_iy_i$ and hence contains all of the generators of $G$. So
$K=G$ contradicting the fact that $K$ is a proper subgroup.
\end{proof}

\section {A non-engulfed proper (locally-free) subgroup}
 In this section we show that every non-elementary word hyperbolic
 group contains subgroups which are not engulfed. More generally we
 show:

\begin{thm} Let $G$ be  a non-elementary word hyperbolic group and $\mathcal F$ a countable collection of quotients of $G$ each with infinite kernel. Then $G$ contains a proper (infinitely generated) subgroup $K$ which surjects on every quotient in the family $\mathcal F$. In particular  $G$ contains a proper subgroup $K$ which is not engulfed.
\end{thm}

\begin{proof}
Enumerate the kernels of the quotients, and for each kernel choose a
set of left coset representatives. Since $G$ is finitely generated
each such set is countable, and we can enumerate the union of the sets
of coset representatives as $g_i, i\in \mathbb N$ with associated
kernels $N_i$.

Choose a proper open subset $\mathcal U$ in $\partial G$.  Since the
kernels are all infinite the limit set of each kernel is dense in the
boundary of $G$. Hence given any finite subset $S_i\subset \mathcal U$
we can choose an infinite order element $y_i\in N_i$ such that
$\partial y_i\subset \mathcal U\setminus S$ and $\partial y_i\cap
\partial g_i=\emptyset$.  Now for sufficiently high powers $y_i^{r_i}$
of $y_i$ and any point $p\not\in \partial y_i$ the image
$y_i^{r_i}g_iy_i^{r_i}(p)$ lies in $\mathcal U$, hence the limit set
of all these elements lies in $\mathcal U$. Setting the subset
$S_i=\bigcup\limits_{j=1}^i \partial y_j$ we may choose these elements
$y_i$ and their powers $r_i$ inductively to ensure that the subset
$\{y_i^{r_i}g_iy_i^{r_i}\mid i=1,\ldots n\}$ freely generates a
subgroup of $G$ with limit set contained in $\mathcal U$, just as we
did in Lemma \ref{torfreegens}. (Again care must be taken over the
choice of absorbing pairs for the elements and we may need to raise
the power of the elements $y_i$.)

It follows that the subgroup generated by any finite subset of these
elements has limit set contained in $\mathcal U$. Any element of the
subgroup $K$ generated by all of these elements lies in one of these
finitely generated subgroups and therefore has its limit set insde
$\mathcal U$. Applying Corollary \ref{limset} we see that $\Lambda K$
is a proper subset of $\partial G$ and so $K$ is a proper (indeed
infinite index) subgroup of $G$.

Consider the image of this subgroup in one of the quotients $G/N\in
\mathcal F$. By construction for each left coset representative $g$ of
the subgroup $N$, the subgroup $K$ contains a generator $y^{r}gy^{r}$
for some element $y\in N$ so $K$ contains a full set of left coset
representatives for each of the kernels in $\mathcal F$ as required.

Now setting $\mathcal F$ to be the set of finite quotients of $G$ we
obtain a proper subgroup which surjects on every finite quotient, and
hence is not engulfed. The ping-pong construction applied at each
stage of the argument shows that we can ensure that the subgroup is an
ascending union of finitely generated free subgroups, and is therefore
locally free.
\end{proof}

Note that the subgroup $K$ constructed in the theorem cannot be
finitely generated since if it were then the ascending chain of
subgroups generated by the finite subsets $\{y_i^{r_i}g_iy_i^{r_i}\mid
i=1,\ldots n\}$ would terminate, which it does not do by construction.

\Addresses\recd
\end{document}